%% file: CMB_final.070712.tex
\documentclass[12pt,dvips]{amsart}
\usepackage{amsfonts, amssymb, latexsym, epsfig}
\usepackage{euler, amsfonts, amssymb, latexsym, epsfig, epic}

\newtheorem{Theorem}{Theorem} 
\newtheorem*{Thm*}{Theorem 1'}

\newtheorem{Proposition}{Proposition} 
\newtheorem{Lemma}{Lemma}

\newtheorem*{Corollary*}{Corollary}
 
\newtheorem*{Theorem*}{Theorem}

\newtheorem*{Question*}{Question}
\theoremstyle{remark}
\newtheorem{Example}{Example}

\theoremstyle{plain}

\newtheorem{Question}{Question}
               
\font\co=lcircle10

\def\jr{\smash{\raise2pt\hbox{\co \rlap{\rlap{\char'005} \char'007}}
               \raise6pt\hbox{\rlap{\vrule height5pt}}
               \raise2pt\hbox{\rlap{\hskip4pt \vrule height0.4pt depth0pt
                width5.7pt}}
               \raise2pt\hbox{\rlap{\hskip-9.5pt \vrule height.4pt depth0pt
                width6.2pt}}
               \lower6pt\hbox{\rlap{\vrule height4.5pt}}}}
\def\je{\smash{\raise2pt\hbox{\co \rlap{\rlap{\char'005}
                \phantom{\char'007}}}\raise6pt\hbox{\rlap{\vrule height5pt}}
               \raise2pt\hbox{\rlap{\hskip-9.5pt \vrule height.4pt depth0pt
                width6.2pt}}}}
\def\er{\smash{\raise2pt\hbox{\co \rlap{\rlap{\phantom{\char'005}} \char'007}}
               \raise2pt\hbox{\rlap{\hskip4pt \vrule height0.4pt depth0pt
                width5.7pt}}
               \lower6pt\hbox{\rlap{\vrule height4.5pt}}}}
\def\+{\smash{\lower6pt\hbox{\rlap{\vrule height17pt}}
                \raise2pt\hbox{\rlap{\hskip-9pt \vrule height.4pt depth0pt
                width18.7pt}}}}
\def\hor{\smash{\raise2pt\hbox{\rlap{\hskip-9.5pt \vrule height.4pt depth0pt
                width19.2pt}}}}
\def\ver{\smash{\lower6pt\hbox{\rlap{\vrule height17pt}}}}

\def\textcross{\ \smash{\lower4pt\hbox{\rlap{\hskip4.15pt\vrule height14pt}}
                \raise2.8pt\hbox{\rlap{\hskip-3pt \vrule height.4pt depth0pt
                width14.7pt}}}\hskip12.7pt}
\def\textelbow{\ \hskip.1pt\smash{\raise2.8pt%
                \hbox{\co \hskip 4.15pt\rlap{\rlap{\char'005} \char'007}
                \lower6.8pt\rlap{\vrule height3.5pt}
                \raise3.6pt\rlap{\vrule height3.5pt}}
                \raise2.8pt\hbox{%
                  \rlap{\hskip-7.15pt \vrule height.4pt depth0pt width3.5pt}%
                  \rlap{\hskip4.05pt \vrule height.4pt depth0pt width3.5pt}}}
                \hskip8.7pt}

\newcommand{\cellsize}{22}
\newlength{\cellsz} 
\newsavebox{\cell}
\sbox{\cell}{\begin{picture}(\cellsize,\cellsize)
\put(0,0){\line(1,0){\cellsize}}
\put(0,0){\line(0,1){\cellsize}}
\put(\cellsize,0){\line(0,1){\cellsize}}
\put(0,\cellsize){\line(1,0){\cellsize}}
\end{picture}}
\newcommand\cellify[1]{\def\thearg{#1}\def\nothing{}%
\ifx\thearg\nothing
\vrule width0pt height\cellsz depth0pt\else
\hbox to 0pt{\usebox{\cell} \hss}\fi%
\vbox to \cellsz{
\vss
\hbox to \cellsz{\hss$#1$\hss}
\vss}}
\newcommand\tableau[1]{\vtop{\let\\\cr
\baselineskip -16000pt \lineskiplimit 16000pt \lineskip 0pt
\ialign{&\cellify{##}\cr#1\crcr}}}
\begin{document}
\pagestyle{plain}

\title{Multiplicity-free Schubert calculus }
\author{Hugh Thomas}
\address{Department of Mathematics and Statistics, University of New Brunswick, Fredericton, New Brunswick, E3B 5A3, Canada }
\email{hugh@math.unb.ca}

\author{Alexander Yong}
\address{Department of Mathematics, University of Minnesota, Minneapolis, MN 55455, USA,\indent{\itshape and}
        The Fields Institute, 222 College Street, Toronto, Ontario,  
        M5T 3J1, Canada}

\email{ayong@math.umn.edu, ayong@fields.utoronto.ca}

\date{July 7, 2007}

%
\maketitle

\section{Introduction}

\subsection{The main result}
Let $Gr(\ell, {\mathbb C}^{n})$ denote the Grassmannian of
$\ell$-dimensional subspaces in ${\mathbb C}^n$. The cohomology ring
${\rm H}^{\star}(Gr(\ell, {\mathbb C}^{n}), {\mathbb Z})$ 
has an additive basis of Schubert classes
$\sigma_{\lambda}$, indexed by Young diagrams $\lambda=(\lambda_1\geq \lambda_2 \geq\ldots \geq \lambda_\ell\geq 0)$ contained in the $\ell\times k$ rectangle
where $k=n-\ell$ (we denote this by $\lambda\subseteq \ell\times k$).
The product of two Schubert classes in 
${\rm H}^{\star}(Gr(\ell, {\mathbb C}^{n}),{\mathbb Z})$
is given by
\begin{equation}
\label{eqn:LR_product}
\sigma_{\lambda}\cdot\sigma_{\mu}=\sum_{\nu\subseteq \ell\times k} c_{\lambda,\mu}^{\nu}\sigma_{\nu},
\end{equation}
where $c_{\lambda,\mu}^{\nu}$ is the classical Littlewood-Richardson 
coefficient (see, e.g.,~\cite{Fulton, ECII}). 

The expansion (\ref{eqn:LR_product}) is {\bf multiplicity-free} if
$c_{\lambda,\mu}^{\nu}\in \{0,1\}$ for all $\nu\subseteq \ell\times k$.
In this paper, we give a nonrecursive, combinatorial answer to 
the following question of W.~Fulton:
\begin{Question*}
When is $\sigma_{\lambda}\cdot\sigma_{\mu}$ multiplicity-free? 
\end{Question*}

Famously, Littlewood-Richardson coefficients also arise as 
decomposition multiplicities of the tensor product
$V^{\lambda}\otimes V^{\mu}=\bigoplus_{\nu}(V^{\nu})^{\oplus c_{\lambda,\mu}^{\nu}}$ of irreducible polynomial representations of $GL(\ell)$.
Earlier, J.~Stembridge~\cite{JRS} solved the analogous question in this 
context; the above question was motivated by this work. (The $GL(\ell)$ problem may 
be regarded as the  special case in the ``$k\to\infty$ limit'' of 
Fulton's question; 
Stembridge's classification is expressed inside our solution below.) 

For partitions 
$\lambda, \mu\subseteq \ell\times k$, place $\lambda$ against the upper left
corner of the rectangle. Then rotate $\mu$ $180$ degrees and place it in the
lower right corner. We refer to ${\rm rotate}(\mu)$ as the resulting
subshape of $\ell\times k$. The {\bf boxes} of these shapes 
refer to their configuration inside $\ell\times k$.

\begin{figure}[htbp]
\begin{center}
\input{configldmu.pstex_t}
\end{center}
\caption{\label{fig:configldmu} $\lambda$ and ${\rm rotate}(\mu)$ inside $\ell\times k$}
\end{figure}
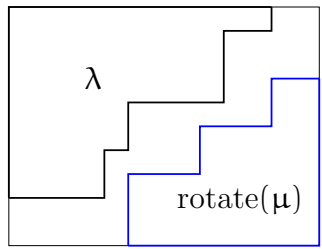

A boring reason for multiplicity-freeness is that  
$\lambda\cap {\rm rotate}(\mu)\neq\emptyset$, since 
the (intersection) product $\sigma_{\lambda}\cdot\sigma_{\mu}$ is merely zero. 
Geometrically, this reflects
the fact that the {\bf Richardson variety} 
$X_{\lambda}^{\mu}=X_{\lambda}(F_{\bullet})\cap X_{\mu}(F_{\bullet}^{{\rm opp}})$ is empty. This variety is the 
scheme-theoretic intersection of Schubert varieties indexed 
by $\lambda$ and $\mu$; these Schubert varieties are in general position
when defined with respect to 
opposite complete flags $F_{\bullet}$ and $F_{\bullet}^{\rm opp}$.

Nevertheless, our main idea is to apply a simple  
extension of this classic 
geometric condition. Define a {\bf Richardson quadruple} to be the
datum $(\lambda,\mu, \ell\times k)$ where $\lambda\cap{\rm rotate}(\mu)=\emptyset$.
If $\lambda\cup {\rm rotate}(\mu)$ does not contain a full 
$\ell$-column or $k$-row, call this Richardson quadruple {\bf basic}.
Otherwise, after removing all full columns and/or
rows, we obtain a Richardson quadruple 
${\widetilde {\mathfrak R}}=({\widetilde \lambda}, {\widetilde \mu}, 
{\widetilde \ell}\times {\widetilde k})$ for smaller partitions
${\widetilde \lambda}, {\widetilde \mu}\subseteq {\widetilde \ell}\times {\widetilde k}$. This latter quadruple is basic, and we call it the 
{\bf basic demolition} of ${\mathfrak R}$.

\begin{Example}
Figure~\ref{fig:configldmu} depicts a basic Richardson quadruple.

The non-basic Richardson quadruple 
${\mathfrak R}=((6,5,4,3,2,1,1), (7,6,6,6,5,2),\, 7\times 9)$
has three full columns ($1,4$ and $5$) and two full rows ($3$ and $4$).
Its basic demolition is 
${\widetilde {\mathfrak R}}=((3,2,1), (5,4,4,2),\, 5\times 6)$:
\begin{figure}[h]
\begin{center}
\epsfig{file=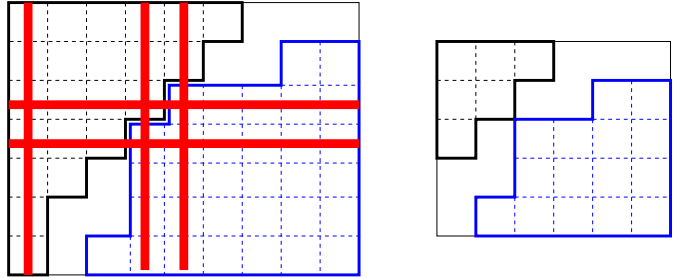, height=1.5in}
\end{center}
\caption{A non-basic Richardson quadruple and its basic demolition}
\end{figure}
\end{Example}
Conceptually, this combinatorics is inspired by a geometric
comparison of Richardson varieties,
$X_{\lambda}^{\mu}\subseteq 
Gr(\ell, {\mathbb C}^{\ell+k})$ and
$X_{\widetilde\lambda}^{\widetilde\mu}\subseteq Gr({\widetilde \ell}, {\mathbb C}^{{\widetilde \ell}
+{\widetilde k}})$. 

In order to state our main result, we need a little 
more notation and terminology, some nonstandard. A {\bf rectangle} is a Young shape
with exactly one distinct part size.
A {\bf fat hook} is a shape whose partition $\lambda$ has two 
distinct part sizes. Furthermore, $\lambda\subseteq \ell\times k$ naturally 
defines a lattice path from the southwest to the northeast corner points 
of the 
rectangle. A {\bf segment} of this lattice path is a maximal consecutive 
sequence of north or east steps. The {\bf shortness} of this
lattice path is the length of its shortest segment. Lastly, a 
{\bf multiplicity-free Richardson quadruple} is a Richardson quadruple
such that the corresponding intersection product $\sigma_{\lambda}\cdot\sigma_{\mu}$ is multiplicity-free.

We are now ready to state our main result:

\begin{Theorem} 
\label{thm:main}
A Richardson quadruple is multiplicity-free if and only if its
basic demolition is multiplicity-free. A basic Richardson quadruple
${\mathfrak R}=(\lambda,\mu,\ell\times k)$ is multiplicity-free if and only if one of the following
conditions holds:
\begin{itemize}
\item[(I)] either $\lambda$ or $\mu$ is a rectangle of shortness 1;
\item[(II)] $\lambda$ is a rectangle of shortness 2, $\mu$ 
is a fat hook (or vice versa);
\item[(III)] $\lambda$ is a rectangle, $\mu$ is a fat hook of shortness 1 (or vice versa);
\item[(IV)] both $\lambda$ and $\mu$ are rectangles.
\end{itemize}
%
%
\end{Theorem}

We remark that once the idea of basic demolitions is found, the solution
to Fulton's question becomes relatively straightforward, 
by exploiting \cite[Theorem~3.1]{JRS}; in the basic case, 
our classification is the same as Stembridge's classification, although
the proof requires further combinatorial analysis. 
However, as is often the case in combinatorics (and we believe here),
the central difficulty of a problem often turns on precisely such
an observation. In fact, after formulating this technique for this work,
we found it to be crucial (in a more general form) for our
arguments in~\cite{HY:comin}.

\begin{Example}
${\mathfrak R}=((4,4,2,2), (3,3,3),\, 6\times 6)$ is \emph{not}
multiplicity-free:
\[\sigma_{(4,4,2,2)}\cdot \sigma_{(3,3,3)} = 2\sigma_{(6,5,4,3,2,1)}+\mbox{multiplicity-free terms},\]
despite the fact that it is a product corresponding to a fat hook
with a rectangle. This is since the lattice path defined by $\lambda$
has shortness~$2$. Meanwhile,
$((4,4,2,2,2), (3,3,3),\, 6\times 6)$ is multiplicity-free, by (III).

Another example is $((4,3,2,1),(4,4,2,2,1),\, 5\times 5)$. This is not basic, but it
is multiplicity-free since its basic demolition is $((1), (1),\, 2\times 2)$.
However $((4,3,2,1), (4,4,2,2,1),\, 6\times 5)$ is not multiplicity-free; 
the tensor product of $GL(6)$ irreducible representations
\[V^{(4,3,2,1)}\otimes V^{(4,4,2,2,1)}=(V^{(5,4,4,4,4,2)})^{\oplus 3}\oplus 
\cdots\]
is not multiplicity-free either. 

Thus, complicated pairs of partitions can give rise
to (nonzero) multiplicity-freeness in the Schubert calculus context. 
This is not true for the $GL(\ell)$-problem (or for basic quadruples): one of the shapes involved
must be a rectangle.
\end{Example}

\subsection{Extensions} 

Investigation of multiplicity-freeness in geometry and representation theory
is of interest, 
see, e.g.,~\cite{Brion, JRS, JRS2} and the references therein. 

In relation to the first of these papers cited, this report 
may be viewed as  part of\linebreak 
``multiplicity-free algebraic geometry''. By this term 
we mean the study
of algebraic subvarieties $Y\subseteq X$ 
that have ``the smallest invariants'' according to the
decomposition of their class into a predetermined 
linear basis of the Chow ring 
$A^{\star}(X)$. (In our case, the class of a
Richardson variety decomposed into Schubert classes for the Grassmannian.)

A feature of our methods is that they extend 
naturally to other (classical) Lie types. In future work, we plan to 
study multiplicity-free Schubert calculus on
cominuscule flag manifolds, a natural generalization 
of ``Grassmannian'', with the following goal:

\begin{Question}
Give a (uniform) characterization of multiplicity-free products of Schubert 
classes, for cominuscule flag manifolds.
\end{Question}
For the main interesting 
cases of Lagrangian and even orthogonal Grassmannians, the
Schubert calculus is determined by the ``shifted tableaux'' combinatorics 
of Schur $P,Q$ polynomials~\cite{JRS:multrule, HB, Pragacz}. 
Briefly, partitions with distinct parts contained in the
staircase $\rho_n=(n, n-1, \ldots, 3,2,1)$ index the Schubert
classes. There is a standard simultaneous placement of shifted
shapes $\lambda$ and $\mu$ into the shifted staircase $\rho_n$. 
When the shapes overlap, the product of their Schubert classes is zero. 
Otherwise, whenever there is a full $i$th column
and $i+1$th row, the situation is a non basic Richardson triple. 
Removal of all such hooks is a basic Richardson triple. 

This allows us to prove multiplicity-free characterizations (omitted here) 
with similar arguments to those found below. Results of C.~Bessenrodt~\cite{Bess}
can replace the role of~\cite{JRS}. (She studies the problem of 
multiplicity-freeness of Schur $P$ polynomials, 
in connection to projective outer products of spin characters.)

Finally, after a version of this paper was made available~\cite{HY}, C.~Gutschwager~\cite{Gut} answered the question
of determining multiplicity-free skew characters. As is explained there,
this problem is equivalent to Fulton's question. 

\section{Demolitions}

In this section, we develop demolition techniques that we will
use in the proof of the main theorem.

\subsection{The Littlewood-Richardson rule} 
We use a standard formulation of
the Littlewood-Richardson rule: $c_{\lambda,\mu}^{\nu}$
counts the number of semistandard fillings of the skew-shape 
$\nu/\lambda$ of content $\mu$ such that the right to left, top 
to bottom reading word $w_1 w_2 \cdots w_{|\mu|}$ 
is a {\bf ballot sequence}, i.e., the number
of appearances of ``$i$'' in $w_1 w_2 \cdots w_j$ is at most the number of 
appearances of
``$i-1$'', for $i\geq 2$ and $1\leq j\leq |\mu|$. 
We call these {\bf LR~fillings}. 

We say a row or column of $\ell\times k$ is {\bf empty} if it neither
contains a box of $\lambda$ nor $\mu$. The following 
{\bf emptiness demolition} is a simple application of the 
Littlewood-Richardson rule:

\begin{Lemma}
\label{lemma:another_redux}
Suppose ${\mathfrak R}$
contains an empty row. Then $\sigma_{\lambda}\cdot\sigma_{\mu}\in {\rm H}^{\star}(Gr(\ell, {\mathbb C}^{\ell+k}))$ is 
the same as $\sigma_{\lambda}\cdot \sigma_{\mu}\in
{\rm H}^{\star}(Gr(\ell-1, {\mathbb C}^{\ell+k-1}))$. In particular,
one product has multiplicity if and only if the other does. A similar
statement holds in the presence of an empty column.
\end{Lemma}
\begin{proof}
This follows the above Littlewood-Richardson rule (and conjugation)
that any $\nu$ that contributes $c_{\lambda,\mu}^{\nu}\neq 0$ satisfies 
$\ell(\nu)\leq \ell(\lambda)+\ell(\mu)\leq \ell-1$ (the latter inequality
being the empty row assumption). So empty rows (or columns)
do not affect the expansion.
\end{proof}

\begin{Example}
The third row in ${\mathfrak R}=((2,1), (2),\, 4\times 3)$ is the only empty 
row/column. The expansion 
\[\sigma_{(2,1)}\cdot \sigma_{(2)} = \sigma_{(2,2,1)}+\sigma_{(3,1,1)}+\sigma_{(3,2)}\]
is the same 
in both ${\rm H}^{\star}(Gr(4, {\mathbb C}^7),\mathbb Z)$ and 
${\rm H}^{\star}(Gr(3, {\mathbb C}^6),\mathbb Z)$.
\end{Example}

\subsection{Basic demolitions} The following lemma explicates the demolition
technique from the statement of the main theorem:

\begin{Lemma}
\label{lemma:basic}
Suppose ${\mathfrak R}$ is not basic, then:
\begin{itemize} 
\item[(a)] ${\widetilde {\mathfrak R}}$ is basic; and
\item[(b)] $\sigma_\lambda \cdot \sigma_\mu$ is multiplicity-free in
${\rm H}^{\star}(Gr(\ell, {\mathbb C}^{\ell+k}))$ if and
only if $\sigma_{\widetilde \lambda} \cdot \sigma_{\widetilde \mu}$ 
is multiplicity-free in
${\rm H}^{\star}(Gr({\widetilde \ell}, {\mathbb C}^{{\widetilde \ell}
+{\widetilde k}}))$.
\end{itemize}
\end{Lemma}
\begin{proof}
For (a), suppose that ${\widetilde {\mathfrak R}}$
contains a full row (the column case is similar). 
Since the boxes in that full row 
were not eliminated by the demolition, the corresponding row in 
${\mathfrak R}$ was not full. Rather these boxes were ``clamped''
together after removing the full columns. But then these boxes,
together with the boxes of the removed full columns from ${\mathfrak R}$,
form a full row in ${\mathfrak R}$, and all would have been removed by
the basic demolition, a contradiction.

For (b),
it is enough to prove the case when ${\widetilde \lambda}$ and
${\widetilde \mu}$ are obtained by removing either one full row or one
full column, since then the stated claim will follow by iterating
this case until the basic situation is reached. 

By conjugation, it suffices
to prove the case that there is a full column in 
$\lambda\cup {\rm rotate}(\mu)$. Now consider the skew shape
$\alpha={\rm rotate}(\mu)^c/\lambda$ where ${\rm rotate}(\mu)^c$ is the complement
in $\ell\times k$ of ${\rm rotate}(\mu)$. Removal of the full column gives
a skew shape ${\widetilde \alpha}$ which has the same number of boxes
as $\alpha$, and also whose boxes are in the same relative position as 
those of $\alpha$. Thus, any LR filling of $\alpha$ of content $\beta$
is also an LR filling of ${\widetilde \alpha}$ of content $\beta$, and
conversely. Moreover, notice that by the Littlewood-Richardson rule, that
in either kind of filling, 
\begin{equation}
\label{eqn:justsmaller}
\beta\subseteq \ell\times (k-1).
\end{equation} 

Finally, we use the well-known symmetry 
\begin{equation}
\label{eqn:symmetry}
c_{\lambda,\mu}^{\nu}=c_{\lambda,{\rm rotate}(\nu)^c}^{{\rm rotate}(\mu)^c}.
\end{equation}
Now assume this number has multiplicity, witnessed by LR fillings of $\alpha$
with content $\beta$, as above. Set $\nu={\rm rotate}(\beta^c)$ 
where rotation and complementation are done
with respect to $\ell\times k$. In view of (\ref{eqn:justsmaller})
we can select ${\widetilde \nu}$ 
such that $\beta={\rm rotate}({\widetilde \nu})^c$, with respect to 
$\ell\times (k-1)$. Now by (\ref{eqn:symmetry}) we see that 
if $c_{\lambda,\mu}^{\nu}\geq 2$ then $c_{{\widetilde\lambda}, {\widetilde\mu}}^{{\widetilde\nu}}\geq 2$ and so if $\sigma_\lambda \cdot \sigma_\mu$ is 
not multiplicity-free then neither is
$\sigma_{\widetilde \lambda} \cdot \sigma_{\widetilde \mu}$. 
The converse argument is similar, starting with LR fillings of 
${\widetilde\alpha}$ of content $\beta$. 
\end{proof}

\subsection{Stembridge demolitions}
We use the following demolition technique throughout the proof
of the main theorem.
Suppose that ${\mathfrak R}$ contains 
a row or column of $\ell\times k$ containing boxes of
$\lambda$, or alternatively ${\rm rotate}(\mu)$ but not both. 
Then define a {\bf Stembridge demolition} to be the Richardson quadruple 
${\overline {\mathfrak R}}=({\overline\lambda}, {\overline\mu}, {\overline \ell}\times {\overline k})$
corresponding to (sequentially) removing such rows or columns. 
(We emphasize that this differs from a basic demolition, which applies 
when a row or column contains \emph{only} boxes of 
(possibly both) $\lambda$ and/or $\mu$.) 

\begin{Example}
${\mathfrak R}=((6,6,4,2), (4,3,2,2),\, 5\times 8)$ is a Richardson
quadruple with multiplicity. 
The columns $1,2,3,4,7$ and $8$
are can be demolished, as can rows $1$ and $5$. If we remove columns $1,2,3$ and row $5$,
we obtain ${\overline {\mathfrak R}}$ having multiplicity. 
However, if we furthermore
remove row $1$, the demolition is multiplicity-free.
\end{Example}

Stembridge demolitions are useful in one direction:

\begin{Lemma}
\label{lemma:stembridge_redux}
If $\sigma_{\overline\lambda}\cdot \sigma_{\overline\mu}\in {\rm H}^{\star}(Gr({\overline \ell}, {\mathbb C}^{{\overline\ell}+{\overline k}}))$ has multiplicity,
then so does $\sigma_{\lambda}\cdot\sigma_{\mu}\in{\rm H}^{\star}(Gr(\ell, {\mathbb C}^{\ell+k}))$.
\end{Lemma}
\begin{proof}
By conjugation, it suffices to handle the case 
that ${\overline \lambda}\cup (r)=\lambda$ (i.e., deleting a row of $\lambda$
gives ${\overline \lambda}$), ${\overline \mu}=\mu$, 
${\overline \ell}=\ell-1$ and ${\overline k}=k$. By assumption, there
exists $\nu\subseteq (\ell-1)\times k$ such that 
$c_{{\overline \lambda},\mu}^{\nu}\geq 2$. 
Moreover, the Littlewood-Richardson rule implies 
$c_{\lambda,\mu}^{\nu\cup (r)}\geq c_{{\overline \lambda},\mu}^{\nu}$: see, 
e.g.,~\cite[Lemma~2.2]{JRS}. Since $\nu\cup (r)\subseteq \ell\times k$, the
claim follows.
\end{proof}

\subsection{A reformulation of Theorem~\ref{thm:main}, inductive Stembridge
demolitions}

We find it useful to give an ``inverse'' formulatation of 
Theorem~\ref{thm:main}, i.e., in terms of when multiplicity
appears:

\begin{Thm*}
\label{thm:grassmannian} 
If ${\mathfrak R}$ is basic then $\sigma_{\lambda}\cdot\sigma_{\mu}$
has multiplicity if and only if:
\begin{itemize}
\item[(I')] $\lambda$ and $\mu$
both have at least two different part sizes; or 

\item[(II')] $\lambda$ has at least three different part sizes and
$\mu$ is a rectangle of shortness at least~2; or

\item[(III')] $\lambda$ is a fat hook of shortness at least~2 and
$\mu$ is a rectangle of shortness at least 3; or

\item[(IV')] cases {\rm ($II^{\prime}$)} or {\rm ($III^{\prime}$)} with the roles of
$\lambda$ and $\mu$ interchanged.
\end{itemize}
Otherwise if ${\mathfrak R}$ is not basic, we conclude as in 
the statement of Theorem~\ref{thm:main}.
\end{Thm*}

This formulation is useful for induction on $\ell,k\geq 1$. For this, we are
interested in the situations where a
Stembridge demolition takes a basic ${\mathfrak R}$
satisfying (I'), (II'), (III') or (IV') of Theorem~1' to
a basic ${\overline {\mathfrak R}}$ that
falls into the same case.
We call these {\bf inductive Stembridge demolitions}.

\begin{Example}
Consider ${\mathfrak R}=((4,4,2,2), (3,3,3),\, 7\times 6)$, which lies in case (III'). 
This has $5$ column and $7$ row Stembridge demolitions available.
However, none are inductive. On the other hand, for ${\mathfrak R}=((4,4,2,2), (4,4,4),\, 7\times 7)$
one can remove columns $5,6$ or $7$ as a inductive Stembridge demolition.
\end{Example}

It is also possible to combine Stembridge and basic demolitions
in inductive arguments, so long as we start and end in the same case
(I')-(IV').

\subsection{Well-ordering corners of fat hooks}
A box of $\lambda$ is a {\bf corner} if there are no
boxes of $\lambda$ below or to the right of it. Similarly, a {\bf corner}
of ${\rm rotate}(\mu)$ is box of this shape without others
above or to the left of it.

Now suppose $\lambda$ and $\mu$ are fat hooks. 
Then both have precisely two corners.
Let $A$ be the lowest/leftmost corner 
of $\lambda$, and $B$ the highest/rightmost.  
Let $X$ and $Y$ be the lowest (equivalently) leftmost corner and
rightmost (equivalently highest) corner of ${\rm rotate}(\mu)$, respectively.  
Let ${\tt row}(A)$ denote the row index of $A$ (as in matrix notation).
Define ${\mathfrak R}$ to be {\bf well-ordered}\footnote{An apology: this
has nothing to do with the usual mathematical notion of being well-ordered} if 
\begin{equation}
\label{eqn:wellordered}
{\tt row}(A)<{\tt row}(X) \mbox{ and } {\tt row}(B)<{\tt row}(Y).
\end{equation} 

\begin{figure}[htb]

\begin{center}
\input{inductstem.pstex_t}
\end{center}
\caption{\label{fig:inductstem} Fat hooks, and their well-ordered corners}
\end{figure}
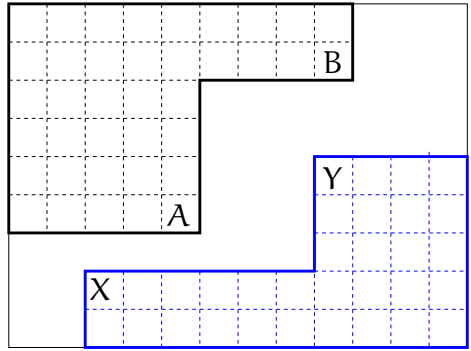

There is an isomorphism
between $Gr(\ell, {\mathbb C}^n)$ and $Gr(k, {\mathbb C}^n)$  
reflected by the operation of {\bf conjugating} the rectangle 
$\ell\times k$ to $k\times \ell$ (and the shapes within). This
sends $\sigma_\lambda$ to $\sigma_{\lambda'}$.
We record the following fact for later use. 

\begin{Proposition}\label{order}
If ${\mathfrak R}$ is basic and $\lambda$ and $\mu$
both have at least two distinct part sizes, then either: 
\begin{itemize}
\item[(a)] ${\mathfrak R}$ has an inductive Stembridge demolition; or
\item[(b)] $\lambda$ or $\mu$ is a hook; or
\item[(c)] $\lambda$ and $\mu$ are fat hooks, and either 
$(\lambda,\mu,\ell\times k)$ or $(\lambda',\mu',k\times \ell)$ is
well-ordered.
\end{itemize}
\end{Proposition}
\noindent
\begin{proof} Basicness implies that
the rightmost column of $\ell\times k$ contains only boxes from 
${\rm rotate}(\mu)$; removing that column is a Stembridge demolition.
However, there are two reasons why this might fail to  
be an {\it inductive} 
Stembridge demolition (for case (I') of Theorem~1'):
\begin{itemize}
\item $\mu$ may no longer have at least two distinct parts; or
\item ${\overline {\mathfrak R}}=(\lambda,{\overline \mu}, \ell\times (k-1))$ 
may not be basic because the top row of $\ell\times (k-1)$  consists 
entirely of boxes from $\lambda$.  
\end{itemize}
(Similar analysis applies to the other three edges 
of $\ell\times k$.)

If $\mu$ has at least three distinct parts then
the first of these possibilities cannot occur.
Thus failure of inductiveness must be
blamed on $\lambda$ extending to the $(k-1)$th column.
Similarly, if removing the bottom row is not inductive,
then $\lambda$ extends all the way to the $(\ell-1)$th row.

Assume further that $\lambda$ is a fat hook. 
If $\lambda$ is not a hook, then the previous paragraph,
together with the basicness implies either removing
the top row or the leftmost column is an inductive Stembridge
demolition. 

Thus, suppose $\lambda$ has at least three distinct 
parts. We may remove the leftmost column of $\ell\times k$
to obtain an inductive Stembridge demolition unless $\mu$ extends
all the way to the second column of $\ell\times k$.  In this case, remove both
the leftmost column and the bottom row.  Since $\mu$ has at least three
parts by assumption, the result of applying these two Stembridge
demolitions is inductive.

We have just disposed of the case where $\mu$ (or symmetrically $\lambda$)
has at least three distinct
parts. Assume now that both $\lambda$ and $\mu$ have 
exactly two distinct parts.  

By (c) we're done if ${\mathfrak R}$ is well-ordered, so 
assume otherwise.  Thus, at least one of 
${\tt row}(X)\leq {\tt row}(A)$ or ${\tt row}(Y)\leq {\tt row}(B)$ holds.  If both of these are true, then the
partition will be well-ordered after conjugating $\ell\times k$.  
If ${\tt row}(X)\leq {\tt row}(A)$ but ${\tt row}(Y)>{\tt row}(B)$  
then it follows that either $\lambda$ is a hook
or there is an inductive Stembridge demolition removing either the top row
or the leftmost column.  Finally, if ${\tt row}(X)>{\tt row}(A)$ but ${\tt row}(Y)\leq {\tt row}(B)$, then
either $\mu$ is a hook or there is an inductive Stembridge demolition obtained 
by removing either the bottom row or the rightmost column.  
\end{proof}

\section{Proof of the main result}

We will actually prove Theorem~1', the reformulation of the main
result given in Section~2. The equivalence of these two statements
is straighforward to check.

In each of the first three subsections below, we assume basicness throughout
and apply the demolitions of
Section~2 to induct on $\ell,k\geq 1$. In the base cases where no induction is 
possible, we describe a skew shape $\nu/\lambda\subseteq \ell\times k$
having two LR~fillings of content $\mu$,
thus showing $c_{\lambda,\mu}^{\nu}\geq 2$. 

\subsection{Proof of multiplicity in case {\rm ($I^{\prime}$)}}
Let $\lambda, \mu$ have at least two distinct nonzero part sizes.
By Proposition~\ref{order}, one of its scenarios (a), (b) or (c)
occurs. We induct if (a) happens. 

Thus consider (b).  Let $\mu=(b,1^a)$ 
be a hook. Let $r$ be the smallest index such that $\lambda_r<\lambda_1$.  
Let $s$ be the smallest index such that $\lambda_s<\lambda_r$. 

For the first filling, 
add a horizontal strip of boxes (i.e., no two in the same column)
of size $b$ with each box labeled ``1'' such that:
\begin{itemize}
\item at least one box appears in the first row, and in the $r$th row; and
\item the maximal possible number of boxes is \emph{not} 
placed in row $s$ (so possibly no boxes occur there) 
\end{itemize}
(Note that by basicness, 
$b<k$ and $\lambda_1'<\ell$,
so such a horizontal strip exists.)
Now add $a$ boxes, no two in the same row, 
with a box at the right-hand end of all rows 
2 through $a+2$, except row $r$.  Label the box in row $i$ with ``$i$''
if $i<r$, and with ``$i-1$'' if $i>r$. 
This clearly gives a LR~filling of a skew 
shape $\nu/\lambda$ with content $\mu$, as desired.

We modify the above to obtain
a second LR~filling:
replace the rightmost ``1'' in row $r$ with an ``$r$''.  
The column below the leftmost ``1'' in row $r$ consists of boxes labeled 
by ``$r$'' to ``$s-2$''.  Increase all of these labels by 1.  
Finally, replace the box labeled ``$s-1$'' in row $s$ by a ``1''.  
\begin{Example}
Being the first of several such arguments, let's take this one in slow-motion.
Consider ${\mathfrak R}=((11,11,11,7,7,4,4,2,2),(12, 1^{9}),\, 11\times 13)$. Here $r=4, s=6$. See Figure~\ref{fig:1} below.
\begin{figure}[h]
\input{caseIb.pstex_t}
\caption{\label{fig:1} Find two LR fillings}
\end{figure}
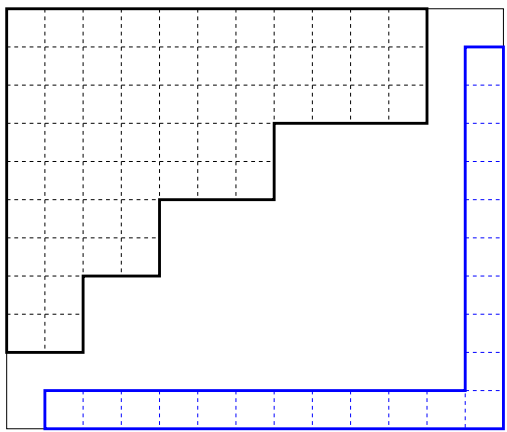
In general there is choice in constructing the two fillings; but not in this 
example. We invite the reader to pencil in the first LR filling, followed
by marking in the modification. The first reading word is
\[1, 1, 2, 3, 1, 1, 1, 1, 4, 5, 1, 1, 6, 7, 1, 1, 8, 9, 1, 1, 10.\] 
The second is 
\[1, 1, 2, 3, 4, 1, 1, 1, 5, 1, 1, 1, 6, 7, 1, 1, 8, 9, 1, 1, 10.\]
How does the filling vary as $a$ and $b$ change?
\end{Example}

Now we turn to (c). Assume without loss of generality that
${\mathfrak R}$ is well-ordered. The example given in
Figure~\ref{fig:inductstem} provides a running example. 

To describe our LR~fillings, 
think of the boxes of ${\rm rotate}(\mu)$ as movable tiles, 
labeled by their row number in $\mu$ (i.e., by $\ell-{\tt row}(\cdot)+1$).
  
Remove the tiles in the columns below $X$ and $Y$ (including
$X$ and $Y$) and set them aside. Shift all the remaining
tiles of ${\rm rotate}(\mu)$ one square to the
left.  Note that since ${\mathfrak R}$ is basic, 
the tiles still don't overlap $\lambda$. 

Using only columns $k-\mu_1 $ to $k-1$, move the $\mu_1 -2$ remaining
columns of tiles (possibly left or right, but maintaining their
relative order) so that
the columns below $A$ and $B$ are empty.\footnote{Depending on the
exact well-ordered configuration of $A,B,X$ and $Y$, there can be choice, but
these all give the same end result. To be precise, we just insist
that the same columns be used for both fillings described.}
Next slide these columns up so they are immediately below (the lattice
path defined by) $\lambda$,
and reverse the order of the tiles in each column. The only worry in this
last step is that there might not be enough room to fit each column.
It is straightforward to check 
that this cannot happen, because of the well-ordered assumption. 

\begin{Example}
\label{exa:8}
Following the running example defined by Figure~\ref{fig:inductstem}, there is
actually no choice in this case 
which columns to use in the ``slide the columns up'' step. The result at that 
stage is depicted below in Figure~\ref{fig:exa8}.

\begin{figure}[h]
\begin{center}
\epsfig{file=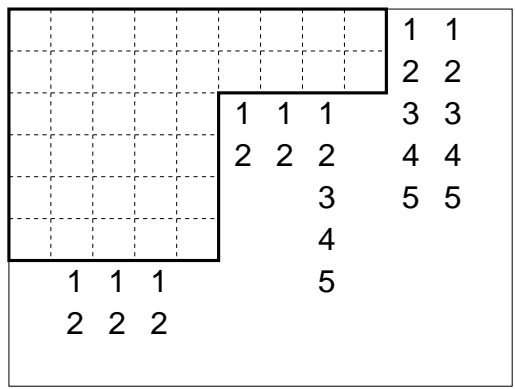, height=2in}
\end{center}
\caption{\label{fig:exa8} In process filling for Example~\ref{exa:8}}
\end{figure}
\end{Example}

Let $a=\ell - {\tt row}(A)$ and $b=\ell -{\tt row}(B)$ count the number of
boxes strictly below $A$ and $B$ respectively. 
Also, let $x=\ell-{\tt row}(X)+1$, $y=\ell-{\tt row}(Y)+1$ 
be the number of boxes in
${\rm rotate}(\mu)$ weakly below $X$ and $Y$. The tiles that we set aside
consist of two copies of $1$ to $x$, together with one copy of $x+1$ to
$y$.  We wish to use these to fill the columns below $A$ and $B$ and the
rightmost column.

For our first filling, put 1 in the rightmost column.  Put 
2 to $y-(a-x)$ in the column below $B$, and put $1$ to $x$ followed
by $y-(a-x)+1$ to $y$ in the column below $A$.  There is space to do this:
under $A$ we use precisely $a$ labels, while under $B$ we use fewer than
$b$ elements, precisely because of well-orderedness (\ref{eqn:wellordered}).

By construction, the associated word is a ballot sequence,
because this is true for the word coming from the columns of labels we've
inserted \emph{not} below $A$ and $B$. (To each ``$i$'', there is an
``$i-1$'' north of it in its column, which appears first in the reading
word.) This property is clearly maintained after adding the
three special 
columns.
The only complaint is that this filling may not yet describe a
skew-shape. We fix this by sliding tiles
to the left (in the same row) as necessary. This doesn't alter the
reading word. Lastly, it is straightforward to verify that the
filling is semistandard. 

\begin{Example}
Continuing our running example, the reader may wish to fill in the extra three
columns. Since $a=3, b=7, x=2, y= 5$, we need to place a column of
``$2, 3, 4$'' under column 9, a column of ``$1,2,5$'' under column 5,
and a single ``1'' in the rightmost column. The final shape after shifting to 
the left is $(12, 11, 11, 11, 9, 7, 5, 4, 1)/\lambda$.
\end{Example}

The second LR filling is similar: put $1$ in the rightmost
column, $1$ to $y-(a-x)-1$ in the column below $B$, and $2$ to 
$x$ followed by $y-(a-x)$ to $y$ in the column below $A$.  
One also checks this satisfies
the necessary conditions, and leads to an LR~filling of 
the same shape as the first filling (even before the final left 
justificiation).  

\begin{Example}
Erasing what we did to columns 9 and 5 in the running example last time, we
put ``$1,2,3$'' under column 9, and ``$2, 4,5$'' under column 5.
\end{Example}

\subsection{Proof of multiplicity in case {\rm ($II^{\prime}$)}}
Let $\lambda$ have at least three part sizes and $\mu$ be a rectangle
of shortness at least two.  Thus, 
$\mu=(g^h)$
where $2\leq g\leq k-2$ and $2\leq h\leq \ell-2$.  

If $g$ or $h$ equals 2;
by conjugation, we may assume that $h=2$. Note that the basicness assumption 
together with the fact that $\lambda$ contains at least three distinct parts, 
implies one can add to $\lambda$ a horizontal strip of length $g+2$ which
uses (at least) four different rows. Call the
rightmost boxes of these four rows $B_1,B_2, B_3$ and $B_4$ 
where ${\tt row}(B_i)<{\tt row}(B_j)$ when $i<j$. 

If $\lambda$ does not extend to the $(\ell-1)$th row, then this horizontal
strip can be done in any way possible; otherwise some extra care has to be
made. In this latter situation, by inductive Stembridge reductions (removing
the first few columns), we may assume that $\lambda$ has at least one
column reaching the $(\ell-1)$ th row. Thus, choose our horizontal strip
to include a (single) box in row $\ell$, and make that box $B_4$.

In either case, label the other $g-2$ boxes with 
``1''s. Place $g-2$ boxes labeled ``2''
under these boxes we have just labeled.
These boxes give a skew shape $\nu/\lambda$ after sliding the boxes
along each row to the left against the lattice path defined by $\lambda$. See
Figure~\ref{fig:caseIIfirst} below.
\begin{figure}[h]
\begin{center}
\epsfig{file=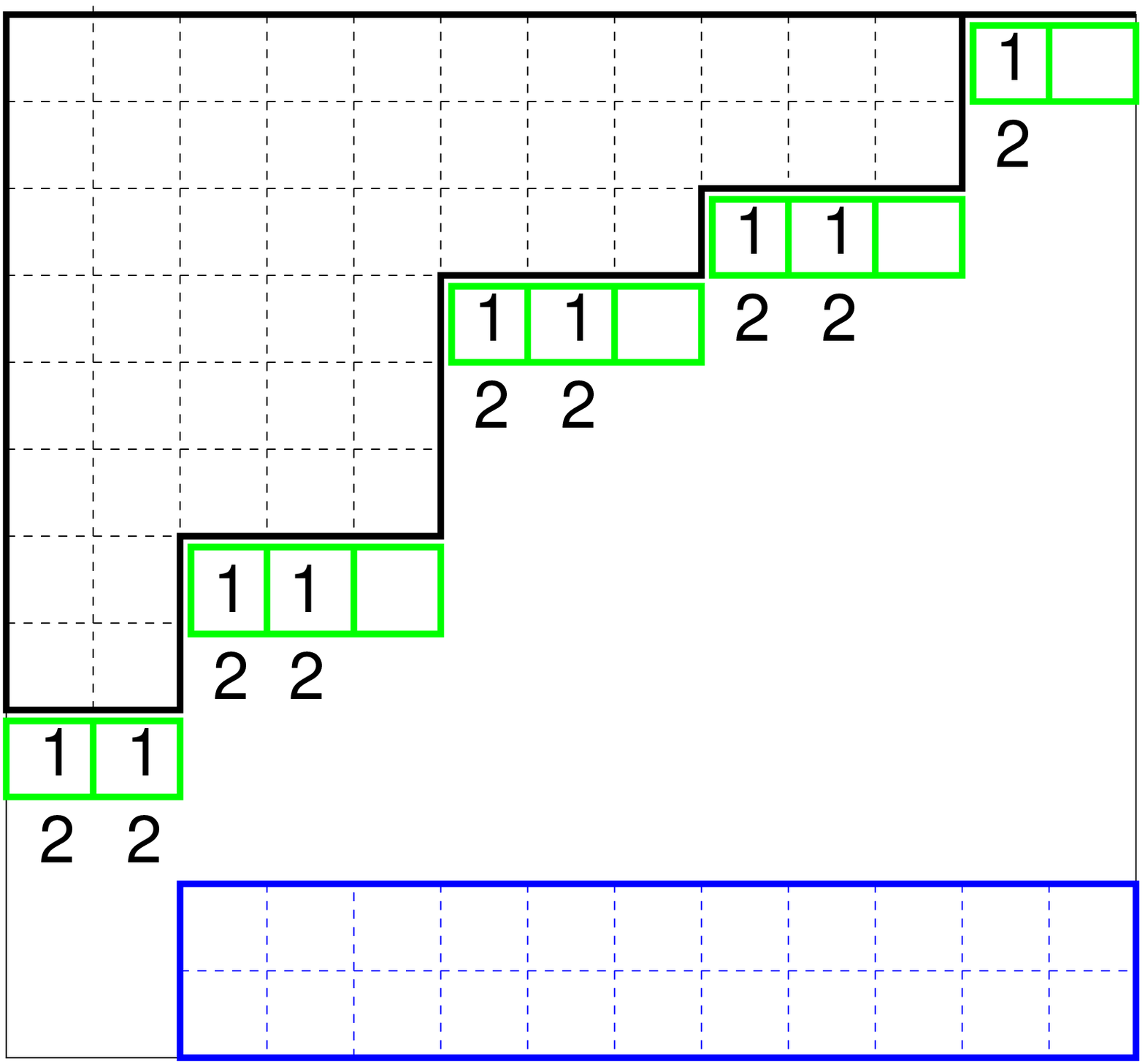, height=2in}
\end{center}
\caption{\label{fig:caseIIfirst} A partial LR filling with blanks $B_1, B_2, B_3, B_4$ to fill}
\end{figure}

One LR~filling of $\nu/\lambda$ is
given by $B_1=B_2=1$, $B_3=B_4=2$; another is given by 
$B_1=B_3=1$, $B_2=B_4=2$. Thus, the result holds if $g$ or $h$ equals~2.

Now assume that $\mu=(g^h)$ where $g,h\geq 3$. 
After applying inductive Stembridge demolitions 
one can assume (after conjugating) that
$\lambda=(k-1,2^a,1^{\ell-a-2})$, for some $a>0$. 
Proof: applying inductive 
Stembridge demolitions to rows and columns of $\lambda$
we can reduce to the case that $\lambda$ has exactly three part sizes. 
If neither the bottom row nor the rightmost column of $\lambda$ can be
removed by an inductive Stembridge reduction, 
$\lambda$ must extend all the way to the $(k-1)$th column
and $(\ell-1)$th row. Further inductive Stembridge demolitions to the first few 
columns/rows of $\lambda$ imply that $\lambda=(k-1, b^a, 1^{\ell-a-2})$. 
Finally, inductive Stembridge demolitions to columns $2$ through $b$,
or rows $2$ through $a$ (depending on the position of $\mu$) allows us to 
deduce $b=2$ or $a=2$. Conjugating we may assume the former.
 
If $a+h+1>\ell$ then the Stembridge demolition obtained by removing the third
column is clearly inductive if $g<k-3$. So we may assume that $g=k-3$. 
But then we can still Stembridge reduce by removing the third column.
The result is non-basic (and thus not inductive). But removing the
full rows that result gives us a smaller case of (II'), so we can still
induct nonetheless. 

\begin{Example}
A minimal example is given by ${\mathfrak R}=((4,2,2,1), (2,2,2),\, 5\times 5)$. 
Removing the third column is a noninductive,
giving ${\overline {\mathfrak R}}=((3,2,2,1), (2,2,2),\, 5\times 4)$ which is not basic,
because of the full third row. The basic demolition results in the
final quadruple ${\widetilde {\overline {\mathfrak R}}}=((3,2,1), (2,2),\, 4\times 4)$, which is in case (II').
\end{Example}

Thus suppose $a+h+1\leq\ell$. We have an
inductive Stembridge demolition by removing the third column
unless $g=k-2$.
Then removing the second row is an inductive Stembridge demolition 
unless $a=1$. 

Summarizing we've reduced to the case that $\lambda=(k-1, 2, 1^{\ell-3})$
and $\mu=((k-2)^h)$. Here is the first LR filling 
of a skew shape $\nu/\lambda$:

\begin{tabular}{ll}
Column $k$:&1\\
Column $k-1$:&$2,\dots,h$\\
Column $i$, $3\leq i \leq k-2$:&$1,\dots, h$\\
Column 2: & $1,\dots,h-1$\\
Column 1: & $h$ 
\end{tabular}

\begin{figure}[h]
\begin{center}
\epsfig{file=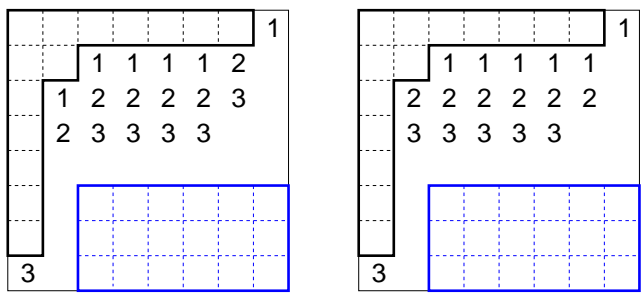, height=1.5in}
\end{center}
\caption{\label{fig:caseIIsecond} Two LR fillings}
\end{figure}

\noindent
For the second filling, interchange columns $2$ and $k-2$. See Figure~\ref{fig:caseIIsecond} for an example.

\subsection{Proof of multiplicity in case {\rm ($III^{\prime}$)}}  
Let $\lambda=(c^d,a^b)$ be a fat hook of shortness at least 2
and let $\mu$ be a rectangle of
shortness at least 3.  Thus, $\mu=(g^h)$
where $3\leq g\leq k-3$ and $3\leq h\leq \ell-3$.

Note that Lemma~\ref{lemma:another_redux} implies that if
$g>k-a$ then we can use emptiness demolitions and basic demolitions
to obtain a basic situation that also
lies in case (III'), so that we can conclude by induction. Similarly, we 
may assume from now on that $h<\ell-d$, since a similar argument holds
if $h>\ell-d$, and since the $h=\ell-d$ reduces to the $g=k-a$ case
handled below, after conjugation.

Next suppose $g=k-a-1$ (we'll deal with $g=k-a$ after this). 
If moreover $3\leq h\leq b$, we have a
filling $F_h$ of a skew shape $\nu/\lambda$
with content $\mu=(g^h)$:

\begin{tabular}{ll}
Column $1$: &$h-1,h$\\
Column $2$: & $h$ \\
Column $c-1$: & $1,\dots,h-2,h$ \\
Column $c$: & $2,\dots,h-1$ \\
Column $k -1$: & $1,\dots,h-1$ \\
Column $k$:  &1\\
All other columns $i>a$:& $1,\dots,h$
\end{tabular}

Our assumptions ensure that the insertion of columns gives rise to a skew-shape
automatically. It is also easy to see that the corresponding reading word
is a ballot sequence. Finally, when checking that the semistandard conditions 
hold, the only worry is between columns $c$ and $c+1$. But notice that by
the assumption that all segments of the lattice path defined by $\lambda$ are
of length at least $2$, the first comparison (in row $d+1>2$) is 
between a ``$2$'' in column $c$ and ``$d+1$'' in column $c+1$
is satisfactory. Since the labels in both columns increment by $1$ as we
go down, all desired inequalities between these two columns are satisfied.
The second filling $G_h$
is obtained by interchanging the bottom-most entries of columns
$c-1$ and $c$.  
It is easy to see that this is also a LR~filling.
In the remaining cases, we invite the reader to consider how the fillings
from the example above are adjusted.
\begin{Example}
Consider ${\mathfrak R}= ((8^2, 3^5), (7^5),\, 12\times 11)$ so $c=8, d=2, a=3, b=5$. Here
$h=5=b$ and the filling given in Figure~\ref{fig:caseIIIfirst}.
\begin{figure}[h]
\begin{center}
\epsfig{file=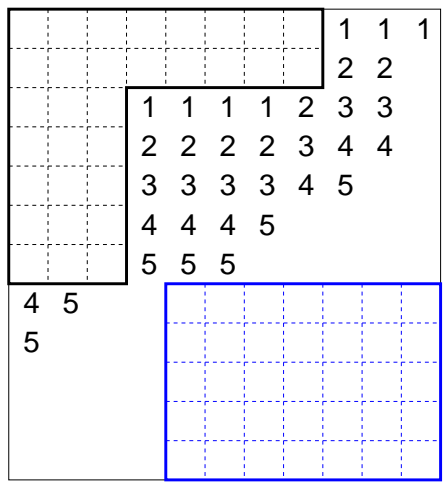, height=2.5in}
\end{center}
\caption{\label{fig:caseIIIfirst} The LR filling $F_5$}
\end{figure}
\end{Example}

If $h=b+1$, then we apply the instructions above, 
but when we place the
entries in the columns as instructed, we see that we do not obtain a skew
shape.  However, we do obtain a skew shape if we push filled boxes to the
left.  Call the resulting fillings $F_{h+1}$ and $G_{h+1}$.  It is easy
to verify that they are again LR~fillings of the same 
shape, of content~$\mu$.  

If $h> b+1$, define $F_h$ and $G_h$ by taking $F_{b+1}$ and $G_{b+1}$ 
and then adding $b+2,\dots,h$ onto the first $g$ columns. (Note that since
$h< \ell-d$, there is always sufficient room to do this.) Now $F_h$ and
$G_h$ clearly serve our purpose.  

Next, when $g=k-a$, note that $a\geq 3$, so we can apply
the same procedures as above, except we move the first two columns
labeled with ``$h-1, h$'' and ``$h$'' to columns 2 and 3, and
insert into column 1 the labels ``$1,\ldots,h$''.
 
Finally, consider what happens when $g<k-a-1$.  If we assume further
that $h\leq b$, then simply remove some columns whose content 
is $1,\dots,h$ from the fillings
$F_h$, $G_h$ defined above, and then move all boxes to the left so as
to form a skew shape.  We apply a similar procedure when $h=b+1$, although
we must remove the columns whose content is $1,\dots,h$
prior to sliding any of the squares over.  When $h>b+1$, build the 
fillings just described for $h=b+1$, and then add $b+2,\dots,h$ to the first
$g$ columns.  

In each case, we produce the two desired LR~fillings.  

\subsection{Conclusion of the proof of Theorems~\ref{thm:main} and~1':}
Having reduced nonbasic quadruples to the basic case, by 
Lemma~\ref{lemma:basic}, we've just shown that the cases (I')-(IV') imply
multiplicity. Conversely, suppose that ${\mathfrak R}$ is basic, but
this quadruple does not satisfy any of (I')-(IV'), then we need to 
show that $\sigma_{\lambda}\cdot\sigma_{\mu}$ is multiplicity-free.

    Since we are not in case (I'), one of $\lambda$ and $\mu$ has
only one part size; so let's assume that $\mu$ is a rectangle.

For a partition 
\[\alpha=(\alpha_1\geq \alpha_2\geq\ldots \geq \alpha_s>0)
\subseteq \ell\times k,\] 
Stembridge \cite{JRS} defines 
\[\alpha^{*}=(\alpha_1-\alpha_\ell \geq \alpha_{1}-\alpha_{\ell-1}\geq \ldots\geq
\alpha_1-\alpha_1).\]

     By conjugation symmetry, we may assume that $\mu$ (as placed in
the bottom right corner of $\ell\times k$) is closer to the top 
edge of $\ell\times k$ than the leftmost edge. We consider the
possibilities for what this distance is. By the basicness assumption, 
it is at least one unit.

     Now, if the distance is precisely one unit, then $\mu^{*}$
is a single row partition, so the
product is multiplicity-free by \cite[Theorem~3.1(i)]{JRS}.

     Next assume that $\mu$ is exactly two units from the top
edge of $\ell\times k$. If $\mu$ is a one-line rectangle (i.e., has only one
nonzero part), the product
is multiplicity-free by \cite[Theorem~3.1(i)]{JRS}, so assume otherwise.  
If $\lambda$ has at least three part sizes,
then we're in case (II'). 
Therefore, $\lambda$ has only one or two part
sizes. Note that in either case, $\lambda$ is a rectangle or fat hook and 
$\mu^{*}$ 
is a two-line rectangle (i.e., it has exactly two equal nonzero part sizes). 
So \cite[Theorem~3.1(ii, iv)]{JRS} allows us
to conclude multiplicity-freeness here also.

     Finally, assume that $\mu$ is at least three units from the top
edge of $\ell\times k$. As before, if $\mu$ is a one-line rectangle, the
product is multiplicity-free, so assume otherwise.  
As above, to avoid being in case (II'),
$\lambda$ has at most two part sizes. If $\mu$ were a two-line
rectangle, the product would be multiplicity-free by 
\cite[Theorem~3.1(ii, iv)]{JRS}, so assume otherwise.
Now avoiding (III') means that 
we
are multiplicity-free, by \cite[Theorem~3.1(iii)]{JRS} ($\lambda$ or, after 
conjugating if necessary, 
$\lambda^{*}$ is a ``near rectangle'' in the terminology 
of~\cite{JRS}).
\qed

\section*{Acknowledgements}
We would like to thank William Fulton and John Stembridge, whose questions
motivated this work. In addition, Stembridge's wonderful symmetric 
functions package ``SF'' (available at 
\textsf{http://www.math.lsa.umich.edu/$\sim$ jrs/maple.html\#SF}) 
was invaluable in the (now hidden) experimental process behind this paper. 

We would also like to thank 
Michel Brion, Calin Chindris, Sergey Fomin,  
Allen Knutson, Ezra Miller, 
Kevin Purbhoo, Vic Reiner, Dennis Stanton, 
Frank Sottile, Terry Tao, Dennis White, Alexander Woo and the
anonymous referees for helpful discussions and suggestions.

AY was partially supported by NSF grant 0601010. 
HT was supported by an NSERC Discovery grant. This work was
partially completed while HT was a visitor, and AY was an NSERC supported visitor, at the Fields Institute 
during the Spring 2005 semester on ``The Geometry of String Theory''; as well
as at the NSF--CBMS Regional Conference on ``Algebraic and
Topological Combinatorics of Ordered Sets'' held
at San Francisco State University in August 2005. 
AY would also like to thank the 2005 AMS Summer Research Institute
on Algebraic Geometry in Seattle, where this work was also carried out.

\end{document}

%% file: configldmu.pstex_t
\begin{picture}(0,0)%
\includegraphics{configldmu.pstex}%
\put(-90,60){$\lambda$}
\put(-55,15){${\rm rotate}(\mu)$}
\end{picture}%
\setlength{\unitlength}{987sp}%
\begingroup\makeatletter\ifx\SetFigFont\undefined%
\gdef\SetFigFont#1#2#3#4#5{%
  \reset@font\fontsize{#1}{#2pt}%
  \fontfamily{#3}\fontseries{#4}\fontshape{#5}%
  \selectfont}%
\fi\endgroup%
\begin{picture}(7834,6034)(1789,-6383)
\end{picture}%

%% file: inductstem.pstex_t
\begin{picture}(0,0)%
\includegraphics{inductstem.pstex}%
\put(-115,47){$A$}
\put(-144,19){$X$}
\put(-56, 61){$Y$}
\put(-56, 105){$B$}
\end{picture}%
\setlength{\unitlength}{1579sp}%
\begingroup\makeatletter\ifx\SetFigFont\undefined%
\gdef\SetFigFont#1#2#3#4#5{%
  \reset@font\fontsize{#1}{#2pt}%
  \fontfamily{#3}\fontseries{#4}\fontshape{#5}%
  \selectfont}%
\fi\endgroup%
\begin{picture}(7288,5488)(2357,-6405)
\end{picture}%

%% file: caseIb.pstex_t
\begin{picture}(0,0)%
\includegraphics{caseIb.pstex}%
\end{picture}%
\setlength{\unitlength}{1579sp}%
\begingroup\makeatletter\ifx\SetFigFont\undefined%
\gdef\SetFigFont#1#2#3#4#5{%
  \reset@font\fontsize{#1}{#2pt}%
  \fontfamily{#3}\fontseries{#4}\fontshape{#5}%
  \selectfont}%
\fi\endgroup%
\begin{picture}(7888,6688)(1757,-7005)
\end{picture}%